\documentclass[11pt]{amsart}

\newtheorem{tw}{Theorem}
\newtheorem{lemma}{Lemma}

\newtheorem{prop}{Proposition}

\newcommand{\cal}{\mathcal}

\title{Twistor spaces of generalized complex structures}
\author{Johann Davidov, Oleg Mushkarov}

\address{Institute of Mathematics and Informatics \\
Bulgarian Academy of Sciences\\ Acad. G.Bonchev Str. Bl.8\\
1113 Sofia\\ Bulgaria} \email{jtd@math.bas.bg, muskarov@math.bas.bg}

\begin{document}

\begin{abstract}
The twistor construction is applied for obtaining examples of generalized complex
structures (in the sense of N. Hitchin) that are not induced by a complex or a
symplectic structure.

 \vspace{0,1cm} \noindent 2000 {\it Mathematics Subject Classification} 53C15, 53C28.

\vspace{0,1cm} \noindent {\it Key words: generalized complex structures, twistor
spaces}

\end{abstract}

\maketitle \vspace{0.5cm}

\section{Introduction}

  The notion of a generalized complex structure has been introduced by N. Hitchin
\cite{Hit}. It generalizes both the concept of a complex structure and that of a
symplectic one and can be considered as the complex analog of the notion of a Dirac
structure introduced by T. Courant and A. Weinstein \cite{Cou,CouWei} to unify Poisson
and symplectic geometries. The generalized complex geometry has been further developed
by M. Gualtieri \cite{Gu} and  has recently attracted the interest of many
mathematicians and physicists, see, for example,
\cite{AbBo,Ba,BaBo,Ber,Ca,Cr,GrMiPeTo-2,Hu,IgWa,Ik,Je,KaLi,Li,LiMiToZa,Zab,Zuc} and the
literature quoted there.

   A generalized complex structure on a smooth manifold $N$ is an endomorphism $J$ of
the bundle $TN\oplus T^{\ast}N$ satisfying the following conditions: $(a)$ $J^2=-Id$,~
$(b)$ $J$ preserves the natural metric $<X+\xi,Y+\eta>=\frac{1}{2}(\xi(Y)+\eta(X))$,
$X,Y\in TN$, $\xi,\eta\in T^{\ast}N$,~ $(c)$ the $+i$-eigensubbundle of $J$ in
$(TN\oplus T^{\ast}N)\otimes {\Bbb C}$ is involutive with respect to the bracket
introduced by T. Courant \cite{Cou}. If $J$ satisfies only the conditions $(a)$ and
$(b)$, it is called generalized almost complex structure. The integrability condition
$(c)$ is equivalent to vanishing of the Nijenhuis tensor of $J$ defined by means of the
Courant bracket instead of the Lie one. Every complex and every symplectic structure on
$N$ induce generalized complex structures on $N$ in a natural way. Examples of
generalized complex structures that cannot be obtained from a complex or a symplectic
structure have been given in \cite{CaGu,Ca,Gu}. The main purpose of the present paper
is to provide other examples of this type by means of the twistor construction.

   The twistor theory has been created by R. Penrose \cite{Pen, PeWa} to solve problems in
Mathematical Physics. The construction of a twistor space in the framework of
Riemannian Geometry has been developed by M. Atiyah, N. Hitchin and I. Singer
\cite{AtHiSi}. Following the Penrose ideas, they have defined a natural almost complex
structure on the twistor space of an oriented  Riemannian four-dimensional manifold and
found its integrability condition. J. Eeells and S. Salamon \cite{EeSa} have introduced
another almost complex structure on the twistor space of a $4$-manifold which although
is never integrable plays an important role in the harmonic maps theory. The twistor
construction has been generalized to any even-dimensional Riemannian manifold by L.
B\'erard-Bergery and T. Ochiai \cite{BeOch}, N. O'Brian and J. Rawnsley \cite{OBrRaw},
P. Dubois-Violette \cite{D-V}, I. Skornyakov \cite{Sk}. It has been extended to the
class of quaternionic K\"ahler manifolds by S. Salamon \cite{Sal}, L. B\'erard-Bergery
(unpublished, see \cite[Sec. 14G]{Besse}, \cite{BeOch}), C. LeBrun \cite{LeBr}.

  Let $M$ be an even-dimensional smooth manifold. Following the general scheme of the
twistor construction, we consider the bundle ${\cal G}$ over $M$ whose fibre at a point
$p\in M$ consists of all generalized complex structures on the vector space $T_pM$
(i.e. endomorphisms satisfying conditions $(a)$ and $(b)$ above) which yield the
canonical orientation of $T_pM\oplus T_p^{\ast}M$. The general fibre of ${\cal G}$
admits a complex structure (in the usual sense) and this fact allows one to define two
natural generalized almost complex structures ${\cal J}_1$ and ${\cal J}_2$ on the
manifold ${\cal G}$ when the base manifold $M$ is endowed with a linear connection.
These are analogs of the Atiyah-Hitchin-Singer and Eells-Salamon almost complex
structures. Suppose that the given connection is torsion-free. Under this condition,
the main result of the present paper states that if $dim\, M=2$, then the structure
${\cal J}_1$ is always integrable, while if $dim\, M\geq 4$, it is integrable if and
only if the given connection on $M$ is flat (i.e. $M$ is an affine manifold). In
contrast, the structure ${\cal J}_2$ is never integrable. The complex structure on the
fibres of ${\cal G}$ is K\" ahlerian with respect to a natural metric induced by the
metric $<~,~>$. The corresponding K\" ahler form yields a generalized complex structure
on the general fibre of ${\cal G}$ and one can define two new generalized almost
complex structures on ${\cal G}$. It is not hard to see that these structures are also
never integrable.

\section{Generalized complex structures}

  Let $W$ be a $n$-dimensional real vector space and $g$ a metric of signature $(p,q)$
on it, $p+q=n$. We shall say that a basis $\{e_1,...,e_{n}\}$ of $W$ is {\it
orthonormal} if $||e_1||^2=...=||e_p||^2=1$, $||e_{p+1}||^2=...=||e_{p+q}||^2=-1$. If
$n=2m$ is an even number and $p=q=m$, the metric $g$ is usually called {\it neutral}.
Recall that a complex structure $J$ on $W$ is called {\it compatible} with the metric
$g$, if the endomorphism $J$ is $g$-skew-symmetric. Suppose that $p=2k$ and $q=2l$, and
let $J$ be a compatible complex structure on $W$. Then it is easy to see by induction
that there is an orthonormal basis $\{e_1,....,e_{n}\}$ of $W$ such that
$e_{2i}=Je_{2i-1}$, $i=1,..,k+l$.

\vspace{0.1cm}

  Now let $V$ be a $2n$-dimensional real vector space and $V^{\ast}$ its dual space. Then
the vector space $W=V\oplus V^{\ast}$ admits a natural neutral metric defined by
\begin{equation}\label{eq 0.0}
<X+\xi,Y+\eta>=\displaystyle{\frac{1}{2}}(\xi(Y)+\eta(X))
\end{equation}
(we refer to \cite{Gu} for algebraic facts about this metric).

\begin{lemma}\label {pi_1,2-basis}
Let $V$ be a $2n$-dimensional real vector space and let $\{e_i+\eta_i\}$, $i=1,...,4n$,
be an orthonormal basis of the space $V\oplus V^{\ast}$ endowed with the neutral metric
{\rm {(\ref{eq 0.0})}}. Then $e_1,...,e_{2n}$ is a basis of $V$ and
$\eta_1,...,\eta_{2n}$ is a basis of  $V^{\ast}$ (similarly for $e_{2n+1},...,e_{4n}$
and $\eta_{2n+1},...,\eta_{4n}$).
\end{lemma}

\begin{proof}
Assume that $\sum_{k=1}^{2n}\lambda_k e_k=0$ for some $\lambda_1,...,\lambda_{2n}\in
{\Bbb R}$. Then we have \\ $\sum_{k=1}^{2n}\lambda_k\eta_l(e_k)=0$ for every
$l=1,...,2n$. Denote by $P$ the matrix $[\eta_l(e_k)]$, $1\leq l,k\leq 2n$. Let $I$ be
the unit $2n\times 2n$-matrix. Then the matrix $S=P-I$ is skew-symmetric since
$\eta_l(e_k)+\eta_k(e_l)=2\delta_{lk}$. Therefore
$(det\,P)^2=det\,(P\,^tP)=det\,(I-S^2)$. Let $\beta_1,...,\beta_{2n}$ be the
eigenvalues of the symmetric matrix $S^2$. Then $\beta_k\leq 0$ since the matrix $S$ is
skew-symmetric. Thus $(det\,P)^2=\Pi_{k=1}^{2n}(1-\beta_k)\geq 1$, in particular
$det\,P\neq 0$. Therefore all $\lambda_k$' s must be zero.

\end{proof}

\vspace{0.1cm}

  Let $\{e_i\}$ be an arbitrary basis  of a real vector space $V$ and $\{\alpha_i\}$ its
dual basis, $i=1,...,2n$. Then the orientation of the space $V\oplus V^{\ast}$
determined by the basis $\{e_i,\alpha_i\}$ does not depend on the choice of the basis
$\{e_i\}$. Further on, we shall always consider $V\oplus V^{\ast}$ with this {\it
canonical orientation}.

 The next lemma is of technical character and will be used in the last section.

\begin{lemma}\label {Orient}
 Let $V$ be a $2$-dimensional real vector space and let $\{Q_i=e_i+\eta_i\}$, $1\leq i\leq 4$,
be an orthonormal basis of the space $V\oplus V^{\ast}$ endowed with its natural
neutral metric {\rm {(\ref{eq 0.0})}}. Then:
$$
\begin{array}{lll}
e_3=a_{11}e_1+a_{12}e_2 \\
e_4=a_{21}e_1+a_{22}e_2
\end{array}
$$
where $A=[a_{kl}]$ is an orthogonal matrix. If $det\,A=1$, the basis $\{Q_i\}$ yields
the canonical orientation of $V\oplus V^{\ast}$ and if $det\,A=-1$ it yields the
opposite one.
\end{lemma}

\begin{proof}
According to Lemma~\ref{pi_1,2-basis},
$$
\begin{array}{lll}
 e_3=a_{11}e_1+a_{12}e_2 \hspace{0.3cm}\mbox { and }\hspace{0.3cm}
                                                                e_1=b_{11}e_3+b_{12}\,e_4\\
                  e_4=a_{21}e_1+a_{22}e_2 \hspace{1.5cm} e_2=b_{21}e_3+b_{22}e_4
\end{array}
$$
where $[b_{kl}]=[a_{kl}]^{-1}$. Since the basis $\{Q_i\}$ is orthonormal, we have
$\eta_i(e_j)+\eta_j(e_i)=0$, $i\neq j$, $1\leq i,j\leq 4$ and
$\eta_1(e_1)=\eta_2(e_2)=1$, $\eta_3(e_3)=\eta_4(e_4)=-1$. The latter identities and
the identities $\eta_1(e_3)+\eta_3(e_1)=0$, $\eta_1(e_4)+\eta_4(e_1)=0$,
$\eta_3(e_4)+\eta_4(e_3)=0$ imply
$$
\begin{array}{c}
a_{12}\eta_1(e_2)+b_{12}\eta_3(e_4)=b_{11}-a_{11} \\
a_{22}\eta_1(e_2)-b_{11}\eta_3(e_4)=b_{12}-a_{21}.
\end{array}
$$
It follows that
\begin{equation}\label{eq I-b-row}
b_{11}^2+b_{12}^2=1.
\end{equation}
Similarly, we see that
\begin{equation}\label{eq rest}
b_{21}^2+b_{22}^2=1,~a_{11}^2+a_{12}^2=1,~a_{21}^2+a_{22}^2=1.
\end{equation}
Expressing $b_{kl}$' s in terms of $a_{kl}$' s, we get from (\ref{eq I-b-row}) and
(\ref{eq rest}) that
$$
(a_{11}a_{22}-a_{12}a_{21})^2=1.
$$
It follows that the matrix $A=[a_{kl}]$ is orthogonal.

  To prove the second part of the lemma, let us denote by $\{\alpha_1,\alpha_2\}$ the dual basis
of the basis $\{e_1,e_2\}$ of $V$. Then
$$
\begin{array}{c}
\eta_1=\alpha_1+c\alpha_2,\> \eta_3=d_{11}\alpha_1+d_{12}\alpha_2\\
\eta_2=-c\alpha_1+\alpha_2,\> \eta_4=d_{21}\alpha_1+d_{22}\alpha_2
\end{array}
$$
for some constants $c$ and $d_{kl}$. For the coefficients $d_{kl}$ we have
$$
\begin{array}{lll}
d_{11}=\eta_3(e_1)=-\eta_1(e_3)=-(a_{11}+ca_{12}),\\
d_{12}=\eta_3(e_2)=-\eta_2(e_3)=-(-ca_{11}+a_{12}),\\
d_{21}=\eta_4(e_1)=-\eta_1(e_4)=-(a_{21}+ca_{22}),\\
d_{22}=\eta_4(e_2)=-\eta_2(e_4)=-(ca_{21}+a_{22}).\\
\end{array}
$$
Thus, if we set $c_{11}=c_{22}=1$ and $c_{12}=-c_{21}=c$, then
$[d_{kl}]=-[a_{kl}][c_{lk}]$. Set $C=[c_{kl}]$ and let $I$ be the unit $2\times
2$-matrix. Then the transition matrix from the basis $\{e_1,e_2,\alpha_1,\alpha_2\}$ to
the basis $\{Q_1,Q_2,Q_3,Q_4\}$ has the form
$$
\left[
\begin{array}{cc}
 I & C \\
 A & -AC^{t}
\end{array} \right]
$$
The determinant of this matrix is equal to $4det\,A$, which proves the lemma.
\end{proof}

\vspace{0.1cm}

  A {\it generalized complex structure} on a real vector space $V$ is, by definition, a
complex structure on the space $V\oplus V^{\ast}$ compatible with its natural neutral
metric \cite{Hit}. If a vector space admits a generalized complex structure, it is
necessarily of even dimension \cite{Gu}. We refer to \cite{Gu} for more facts about the
generalized complex structures.

\vspace{0.1cm}

\noindent {\bf Examples}~\cite{Gu,Hit}. Let $V$ be a $2n$-dimensional real vector
space.

\noindent{\bf 1}. Let $K$ be a complex structure on $V$ and define a complex structure
$K^{\ast}$ on $V^{\star}$ by setting $(K^{\ast}\alpha)(X)=\alpha(KX)$, $\alpha\in
V^{\ast}$, $X\in V$. Then the endomorphism $J$ on $V\oplus V^{\ast}$ defined by $J=K$
on $V$ and $J=-K^{\ast}$ on $V^{\star}$ is a generalized complex structure on $V$. This
structure  yields the canonical orientation of $V\oplus V^{\ast}$ (the orientation
induced by $J$ is defined by means of a basis of the form
$Q_1,JQ_2,...,Q_{2n},JQ_{2n}$) .

\noindent{\bf 2}. Let $\omega$ be a symplectic form on $V$ (i.e. a non-degenerate
$2$-form). Then the map $X\to \imath_X\omega$ is an isomorphism of $V$ onto
$V^{\star}$. Denote this isomorphism also by $\omega$ and define a complex structure
$S$ on $V\oplus V^{\ast}$ by setting $SX=\omega(X)$ and $S\alpha=-\omega^{-1}(\alpha)$
for $X\in V$ and $\alpha\in V^{\ast}$. Then $S$ is compatible with the natural neutral
metric of $V\oplus V^{\ast}$, so $S$ is a generalized complex structure on $V$. The
structure $S$ induces the canonical orientation of $V\oplus V^{\ast}$ if and only if
$n=\frac{1}{2}dim\,V$ is an even number.

  Now let $g$ be a metric on $V$ (of any signature) and $K$ a complex structure on $V$
compatible with the metric $g$. Then the generalized complex structures $J$ and $S$
yielded by $K$ and the $2$-form $\omega(X,Y)=g(KX,Y)$ commute.

\noindent{\bf 3}. The direct sum of generalized complex structures is also a
generalized complex structure.

\noindent{\bf 4}. Any $2$-form $B\in \Lambda^{2}V^{\ast}$ acts on $V\oplus V^{\ast}$
via the inclusion $\Lambda^{2}V^{\ast}\subset \Lambda^{2}(V\oplus V^{\ast})\cong
so(V\oplus V^{\ast})$; in fact this is the action $X+\xi\to \imath_{X}B$;~ $X\in V$,
$\xi\in V^{\ast}$. Denote the latter map again by $B$. Then the invertible map $e^{B}$
is given by $X+\xi\to X+\xi+\imath_{X}B$ and is an orthogonal transformation of
$V\oplus V^{\ast}$. Thus, given a generalized complex structure $J$ on $V$, the map
$e^{B}Je^{-B}$ is also a generalized complex structure on $V$, called the $B$-transform
of $J$.

   Similarly, any $2$-vector $\beta\in \Lambda^{2}V$ acts on $V\oplus V^{\ast}$. If we
identify $V$ with $(V^{\ast})^{\ast}$, so $\Lambda^{2}V\cong
\Lambda^{2}(V^{\ast})^{\ast}$, the action is given by $X+\xi\to \imath_{\xi}\beta\in
V$. Denote this map by $\beta$. Then the exponential map $e^{\beta}$ acts on $V\oplus
V^{\ast}$ via $X+\xi\to X+\imath_{\xi}\beta+\xi$, in particular $e^{\beta}$ is an
orthogonal transformation. Hence, if $J$ is a generalized complex structure on $V$, so
is $e^{\beta}Je^{-\beta}$. It is called the $\beta$-transform of $J$.

\vspace{0.1cm}

   Let $W$ be a $2m$-dimensional real vector space equipped with a metric $g$ of
signature $(2p,2q)$, $p+q=m$. Denote by $J(W)$ the set of all complex structures on $W$
compatible with the metric $g$. The group $O(g)$ of orthogonal transformations of $W$
acts transitively on $J(W)$ by conjugation and $J(W)$ can be identified with the
homogeneous space $O(2p,2q)/U(p,q)$. In particular, ${\it dim}\,J(W)=m^2-m$. The group
$O(2p,2q)$ has four connected components, while $U(p,q)$ is connected, therefore $J(W)$
has four components.

\vspace{0.1cm}

\noindent {\bf Example 5}. The space $O(2,2)/U(1,1)$ is the disjoint union of two
copies of the hyperboloid $x_1^2-x_2^2-x_3^2=1$. It seems instructive to see this in
the context of compatible complex structures. Let $W$ be a $4$-dimensional real vector
space equipped with a neutral metric and $e_1,...,e_4$ an orthonormal basis of $W$. Set
$\varepsilon_k=||e_i||^2$, $k=1,...,4$, and define skew-symmetric endomorphisms of $W$
by setting $S_{ij}e_k=\varepsilon_k(\delta_{ik}e_j - \delta_{kj}e_i)$, $1\leq i,j,k\leq
4$. Then the endomorphisms
$$
\begin{array}{lll}
I_1=S_{12}-S_{34}, & & J_1=S_{12}+S_{34},  \\
I_2=S_{13}-S_{24}, & & J_2=S_{13}+S_{24},  \\
I_3=S_{14}+S_{23}, & & J_3=S_{14}-S_{23}
\end{array}
$$
constitute a basis of the space of skew-symmetric endomorphisms of $W$ subject to
following relations: $I_1^2=-Id$, $I_2^2=I_3^2=Id$, $J_1^2=-Id$, $J_2^2=J_3^2=Id$,
$I_rI_s=-I_sI_r$, $J_rJ_s=-J_sJ_r$, $1\leq r\neq s\leq 3$ and $I_rJ_s=I_sJ_r$, $1\leq
r,s\leq 3$. Let $K$ be a complex structure on W compatible with the metric and let us
set $K=\sum_{r=1}^{3}(x_rI_r+y_rJ_r)$.  Then we have
$$
(-x_1^2+x_2^2+x_3^2-y_1^2+y_2^2+y_3^2)Id+2\sum_{r,s=1}^{3}x_ry_sI_rJ_s=-Id
$$
Evaluating the latter identity at $e_1,...,e_4$, we see that
$$
-x_1^2+x_2^2+x_3^2-y_1^2+y_2^2+y_3^2=-1 \mbox { and } x_ry_s=0 \mbox { for } r,s=1,2,3.
$$
Therefore $K^2=-Id$ if and only if either $x_1^2-x_2^2-x_3^2=1$ and $y_1=y_2=y_3=0$ or
$y_1^2-y_2^2-y_3^2=1$ and $x_1=x_2=x_3=0$.

\vspace{0.1cm}

  Consider $J(W)$ as a (closed) submanifold of the vector space $so(g)$ of $g$-skew-symmetric
endomorphisms of $W$. Then the tangent space of $J(W)$ at a point $J$ consists of all
endomorphisms $Q\in so(g)$ anti-commuting with $J$. Thus we have a natural $O(g)$ -
invariant almost complex structure ${\cal K}$ on $J(W)$ defined by ${\cal K}Q=J\circ
Q$. It is easy to check that this structure is integrable.

 Fix an orientation on $W$ and denote by $J^{\pm}(W)$ the set of compatible complex
structures on $W$ that induce $\pm$ the orientation of $W$. The set $J^{\pm}(W)$ has
the homogeneous representation $SO(2p,2q)/U(p,q)$ and, thus, is the union of two
components of $J(W)$.

\vspace{0.1cm}

\noindent {\bf Example 6}. Under the notations of Example 5, let $e_1,e_2,e_3,e_4$ be
an oriented orthonormal basis of $W$. Then it is easy to see that $J^{+}(W)$ is the
hyperboloid $\{\sum_{r=1}^{3}x_rI_r:~x_1^2-x_2^2-x_3^2=1\}$.

\vspace{0.1cm}

 Further on, given an even-dimensional real vector space $V$, the set $J^{+}(V\oplus
V^{\ast})$ of generalized complex structures on $V$ inducing the canonical orientation
of $V\oplus V^{\ast}$ will be denoted by $G(V)$.

 The group $GL(V)$ acts on $V\oplus V^{\ast}$ by letting $GL(V)$ act on $V^{\ast}$ in
the standard way. This action preserves the neutral metric (\ref{eq 0.0}) and the
canonical orientation of $V\oplus V^{\ast}$. Thus, we have an embedding of $GL(V)$ into
the group $SO(<~,~>)$ and, via this embedding, $GL(V)$ acts on the manifold $G(V)$ in a
natural manner.

\vspace{0.1cm}

{\it A generalized almost complex structure} on an even-dimensional smooth manifold $N$
is, by definition, an endomorphism $J$ of the bundle $TN\oplus T^{\ast}N$ with
$J^2=-Id$ which preserves the natural neutral metric of $TN\oplus T^{\ast}N$. Such a
structure is said to be {\it integrable} or {\it a generalized complex structure} if
its $+i$-egensubbunle of $(TN\oplus T^{\ast}N)\oplus {\Bbb C}$ is closed under the
Courant bracket \cite{Hit}. Recall that if $X,Y$ are vector fields on $N$ and
$\xi,\eta$ -- $1$-forms, the Courant bracket \cite{Cou} is defined by the formula
$$
[X+\xi,Y+\eta]=[X,Y]+{\cal L}_{X}\eta-{\cal
L}_{Y}\xi-\frac{1}{2}d(\imath_X\eta-\imath_Y\xi),
$$
where $[X,Y]$ on the right hand-side is the Lie bracket and ${\cal L}$ means the Lie
derivative. As in the case of almost complex structures, the integrability condition
for a generalized almost complex structure $J$ is equivalent to the vanishing of its
Nijenhuis tensor $N$, the latter being defined by means of the Courant bracket:
$$
N(A,B)=-[A,B]-J[A,JB]-J[JA,B]+[JA,JB], \> A,B\in TN\oplus T^{\ast}N.
$$

\vspace{0.1cm}

\noindent {\bf Example 7}. Let $J$ be a generalized almost complex structure on a
manifold $N$ and let $B$ be a smooth $2$-form on $N$. Then, according to Example 4,
$e^{B}Je^{-B}$ is a generalized almost complex structure on $N$. The exponential map
$e^{B}$ is an authomorphism of the Courant bracket (i.e.
$[e^{B}(X+\xi),e^{B}(Y+\eta]=e^{B}[X+\xi,Y+\eta]$) if and only if the form $B$ is
closed \cite{Gu}.  In this case the structure $e^{B}Je^{-B}$ is integrable exactly when
the structure $J$ is so.

\section{The twistor space of generalized complex structures}

  Let $M$ be a smooth manifold of dimension $2n$. Denote
by $\pi:{\cal G}\to M$ the bundle over $M$ whose fibre at a point $p\in M$ consists of
all generalized complex structures on $T_pM$ that induce the canonical orientation of
$T_pM\oplus T^{\ast}_pM$. This is the associated bundle
$$
GL(M)\times_{GL(2n,{\Bbb R}
)} G({\Bbb R}^{2n})
$$
where $GL(M)$ denotes the principal bundle of linear frames on $M$.

  Let $\nabla$ be a linear connection on $M$. Following the standard twistor
construction,  we can define two generalized almost complex structures ${\cal
J}_{1}^{\nabla}$ and ${\cal J}_{2}^{\nabla}$ on the manifold ${\cal G}$ in the
following way: The connection $\nabla$ gives rise to a splitting ${\cal V}\oplus {\cal
H}$ of the tangent bundle of the associated bundle ${\cal G}$ into vertical and
horizontal parts. The vertical space ${\cal V}_J$ of ${\cal G}$ at a point $J\in {\cal
G}$ is the tangent space at $J$ of the fibre through this point. This fibre is the
manifold $G(T_pM)$, $p=\pi(J)$, which admits a natural complex structure ${\cal K}$
defined in the previous section and we set ${\cal
J}_{\alpha}^{\nabla}=(-1)^{\alpha+1}{\cal K}$ on ${\cal V}_J$ and  ${\cal
J}_{\alpha}^{\nabla}=(-1)^{\alpha}{\cal K}^{\ast}$ on ${\cal V}_J^{\ast}$,
$\alpha=1,2$. Thus ${\cal J}_{\alpha}^{\nabla}U=(-1)^{\alpha+1}J\circ U$ for every
$U\in {\cal V}_J$ ($U$ being considered as an endomorphism of $T_pM\oplus T^{\ast}_pM$)
and $({\cal J}_{\alpha}^{\nabla}\omega)(U)=(-1)^{\alpha}\omega(J\circ U)$ for
$\omega\in {\cal V}_J^{\ast}$. The horizontal space ${\cal H}_J$ is isomorphic via the
differential $\pi_{\ast J}$ to the tangent space $T_{p}M, p=\pi(J)$. Denoting
$\pi_{\ast J}|{\cal H}$ by $\pi_{{\cal H}}$, we define ${\cal J}_{\alpha}^{\nabla}$ on
${\cal H}_J\oplus {\cal H}_J^{\ast}$ as the lift of the endomorphism $J$ by the map
$\pi_{{\cal H}}\oplus (\pi_{{\cal H}}^{-1})^{\ast}$.

\vspace{0.1cm}

\noindent

{\bf Remark}. The fibre of ${\cal G}$ at any point $p\in M$ contains generalized
complex structures on $T_pM$ which do not preserve $T_pM$ as well as structures which
do not send $T_pM$ onto $T_{p}^{\ast}M$. This shows that the generalized almost complex
structures structures ${\cal J}_{1}^{\nabla}$ and ${\cal J}_{2}^{\nabla}$ are not
induced by an almost complex or a symplectic structure.

\vspace{0.1cm}

    Further the generalized almost complex structure ${\cal J}_{\alpha}^{\nabla}$ will be
simply denoted by ${\cal J}_{\alpha}$ when the connection $\nabla$ is understood. The
image of every $A\in T_pM\oplus T_p^{\ast}M$ under the map $\pi_{{\cal H}}^{-1}\oplus
\pi_{{\cal H}}^{\ast}$ will be denoted by $A^h$. The elements of ${\cal H}_J^{\ast}$,
resp. ${\cal V}_J^{\ast}$, will be considered as $1$-forms on $T_J{\cal G}$ vanishing
on ${\cal V}_J$, resp. ${\cal H}_J$.

\vspace{0.1cm}

 Let $A(M)$ be the bundle of the endomorphisms of $TM\oplus T^{\ast}M$ which are skew-
symmetric with respect to its natural neutral metric $<~,~>$.  Consider the twistor
space ${\cal G}$ as a subbundle of $A(M)$. Then the inclusion of ${\cal G}$ is
fibre-preserving and the horizontal space of ${\cal G}$ at a point $J$ coincides with
the horizontal space of $A(M)$ at that point since the inclusion $G({\Bbb
R}^{2n})\subset so(2n,2n)$ is $SO(2n,2n)$-equivariant. Let $(U,x_1,...,x_{2n})$ be a
local coordinate system of $M$ and $\{Q_1,...,Q_{4n}\}$ an orthonormal frame of
$TM\oplus T^{\ast}M$ on $U$. Set $\varepsilon_k=||Q_k||^2$, $k=1,...,4n$, and define
sections $S_{ij}$, $1\leq i,j\leq {4n}$, of $A(M)$ by the formula
\begin{equation}\label{eq Sij}
S_{ij}Q_k=\varepsilon_k(\delta_{ik}Q_j - \delta_{kj}Q_i).
\end{equation}
Then $S_{ij}$, $i<j$, form an orthogonal frame of $A(M)$ with respect to the metric
$<a,b>=\displaystyle{-\frac{1}{2}}Trace\,(a\circ b);\, a,b\in A(M)$; moreover
$||S_{ij}||^2=\varepsilon_i\varepsilon_j$ for $i\neq j$.  Set
$$\tilde x_{l}(a)=x_{l}\circ\pi(a),~ y_{ij}(a)=\varepsilon_i\varepsilon_j<a,S_{ij}>$$ for
$a\in A(M)$. Then $(\tilde x_{l},y_{ij})$, $1\leq l\leq 2n$, $1\leq i < j\leq 4n$, is a
local coordinate system of the manifold $A(M)$.

Let
$$
V=\sum_{i<j}v_{ij}\frac{\partial}{\partial y_{ij}}(J)
$$
be a vertical vector of ${\cal G}$ at a point $J$.  It is convenient to set
$v_{ij}=-v_{ji}$ for $i\geq j$, $1\leq i,j\leq {4n}$. Then the endomorphism $V$ of
$T_{p}M\oplus T_{p}^{\ast}M$, $p=\pi(J)$, is determined by $VQ_i=\sum_{j=1}^{4n}
\varepsilon_iv_{ij}Q_j$. Moreover
\begin{equation}\label{cal J/ver}
{\cal J}_{\alpha}V=(-1)^{\alpha+1}\sum_{i<j}
(\sum_{k=1}^{4n}v_{ik}y_{kj}\varepsilon_k)\frac{\partial}{\partial y_{ij}}.
\end{equation}

 Note also that, for every $A\in T_{p}M\oplus T_{p}^{\ast}M$, we have

\begin{equation}\label{cal J/hor}
A^h=\sum_{i=1}^{4n}(<A,Q_i>\circ\pi)\varepsilon_i Q_i^h  \mbox {  and  } {\cal
J}_{\alpha}A^h=\sum_{i,j=1}^{4n}(<A,Q_i>\circ\pi)y_{ij}Q_j^h.
\end{equation}

  For each vector field
$$X=\sum_{i=1}^{2n} X^{i}\frac{\partial}{\partial x_i}$$
on $U$, the horizontal lift $X^h$ on $\pi^{-1}(U)$ is given by
\begin{equation}\label{eq 3.1}
X^{h}=\sum_{k=1}^{2n} (X^{l}\circ\pi)\frac{\partial}{\partial\tilde x_l}-
\sum_{i<j}\sum_{p<q}
y_{pq}(<\nabla_{X}S_{pq},S_{ij}>\circ\pi)\varepsilon_i\varepsilon_j\frac{\partial}{\partial
y_{ij}}
\end{equation}
where $\nabla$ is the connection on $A(M)$ induced by the given connection on $M$.

     Let $a\in A(M)$ and $p=\pi(a)$. Then (\ref{eq 3.1}) implies
that, under the standard identification of $T_{a}A(M)$ with the vector space of
skew-symmetric endomorphisms of $T_{p}M\oplus T_{p}^{\ast}M$, we have
\begin{equation}\label{eq 3.2}
[X^{h},Y^{h}]_{a}=[X,Y]^h_a + R(X,Y)a
\end{equation}
where $R(X,Y)a$ is the curvature of the connection $\nabla$ on $A(M)$ (for the
curvature tensor we adopt the following definition: $R(X,Y)=\nabla_{[X,Y]}-
[\nabla_{X},\nabla_{Y}]$).

{\it Notations}. Let $J\in {\cal G}$ and $p=\pi(J)$. Take an oriented orthonormal basis
$\{a_1,...,a_{4n}\}$ of $T_pM\oplus T_{p}^{\ast}M$ such that $a_{2l}=Ja_{2l-1}$,
$l=1,..,2n$. Let $\{Q_i\}$, $i=1,...,4n$, be an oriented orthonormal frame of $TM\oplus
T^{\ast}M$ near the point $p$ such that
$$
Q_i(p)=a_i \mbox { and } \nabla Q_i|_p=0, \> i=1,...,4n.
$$

  Define a section $s$ of ${\cal G}$ by setting
$$
sQ_{2l-1}=Q_{2l}, ~ sQ_{2l}=-Q_{2l-1},\> l=1,...,2n.
$$
Then, considering $s$ as a section of $A(M)$, we have
$$
s(p)=J, ~ \nabla s|_p=0.
$$
Thus $X^h_J=s_{\ast}X$ for every $X\in T_pM$.

 Further, given a smooth manifold $N$, the natural projections of $TN\oplus T^{\ast}N$
onto $TN$ and $T^{\ast}N$ will be denoted by $\pi_1$ and $\pi_2$, respectively.

  We shall use the above notations throughout this section.

\vspace{0.1cm}

To compute the Nijenhuis tensor of the generalize almost complex structure ${\cal
J}_{\alpha}$, $\alpha=1,2$, on the twistor space ${\cal G}$ we need some preliminary
lemmas.

\vspace{0.1cm}

\begin{lemma}\label {brackets}
If $A$ and $B$ are sections of the bundle $TM\oplus T^{\ast}M$ near $p$, then:
\begin{enumerate}
\item[$(i)$]
$[\pi_1(A^h),\pi_1({\cal
J}_{\alpha}B^h)]_J=[\pi_1(A),\pi_1(sB)]^h_J+R(\pi_1(A),\pi_1(JB))J.$ \vspace{0.2cm}
\item[$(ii)$]
$[\pi_1({\cal J}_{\alpha}A^h),\pi_1({\cal
J}_{\alpha}B^h)]_J=[\pi_1(sA),\pi_1(sB)]^h_J+R(\pi_1(JA),\pi_1(JB))J.$
\end{enumerate}
\end{lemma}

\begin{proof}
Set $X=\pi_1(A)$. By (\ref{eq 3.1}), we have  $X_J^h=\sum_{l=1}^{2n}
X^l(p)\frac{\partial}{\partial\tilde x_l}(J)$ since  $\nabla S_{ij}|_p=0$,
$i,j=1,...,4n$. Then, using (\ref{cal J/hor}), we get
$$
[X^h,\pi_1({\cal J}_{\alpha}B^h)]_J=
$$
$$
\sum_{i,j=1}^{4n}<B,Q_i>_py_{ij}(J)[X^h,\pi_1(Q_j)^h]_J +
X_p(<B,Q_i>)y_{ij}(J)(\pi_1(Q_j))^h_J.
$$
Moreover, $sB=\sum_{ij}<B,Q_i>(y_{ij}\circ s)Q_j$ since $({\cal J}_{\alpha}B^h)\circ
s=(sB)^h\circ s$. Now formula $(i)$ follows by means of (\ref{eq 3.2}). Similar
computations give $(ii)$.
\end{proof}

\begin{lemma}\label {Lie deriv}
Let $A$ and $B$ be sections of the bundle $TM\oplus T^{\ast}M$ near $p$, and let $Z\in
T_pM$, $W\in {\cal V}_J$. Then:
\begin{enumerate}
\item[$(i)$]
\hspace{0.2cm}$({\cal L}_{\pi_1(A^h)}{\pi_2(B^h)})_J=({\cal
L}_{\pi_1(A)}{\pi_2(B)})^h_J.$ \vspace{0.2cm}
\item[$(ii)$]
\hspace{0.2cm}$({\cal L}_{\pi_1(A^h)}{\pi_2({\cal J}_{\alpha}B^h)})_J=({\cal
L}_{\pi_1(A)}{\pi_2(sB)})^h_J.$ \vspace{0.2cm}
\item[$(iii)$] \ \\
$
\begin{array}{lll}
({\cal L}_{\pi_1({\cal J}_{\alpha}A^h)}\pi_2(B^h))_J(Z^h+W)=\\
                                                    \\
({\cal L}_{\pi_1(sA)}\pi_2(B))^h_J(Z^h)+(\pi_2(B))_p(\pi_1(WA)).
\end{array}
$
 \vspace{0.2cm}
\item[$(iv)$] \ \\
$
\begin{array}{lll}
({\cal L}_{\pi_1({\cal J}_{\alpha}A^h)}\pi_2({\cal J}_{\alpha}B^h))_J(Z^h+W)=\\
                                                             \\
({\cal L}_{\pi_1(sA)}\pi_2(sB))^h_J(Z^h)+(\pi_2(JB))_p(\pi_1(WA)).
\end{array}
$
\end{enumerate}
\end{lemma}

\begin{proof}
Formula $(i)$  follows from (\ref{eq 3.2}). Simple computations involving (\ref{cal
J/hor}) and (\ref{eq 3.2}) give $(ii)$, $(iii)$ and $(iv)$.
\end{proof}

\vspace{0.1cm}

  The next lemma is also easy to prove by means of (\ref{cal J/hor}) and (\ref{eq
  3.2}).

\begin{lemma}\label {Half-diff}
Let $A$ and $B$ are sections of the bundle $TM\oplus T^{\ast}M$ near $p$. Let $Z\in
T_pM$ and $W\in {\cal V}_J$. Then:
\begin{enumerate}
\item[$(i)$]
\hspace{0.2cm}$(d\>\imath_{\pi_1(A^h)}\pi_2(B^h))_J=(d\>\imath_{\pi_1(A)}\pi_2(B))^h_J$
\vspace{0.2cm}
\item[$(ii)$] \ \\
$
\begin{array}{lll}
(d\>\imath_{\pi_1(A^h)}\pi_2({\cal J}_{\alpha} B^h))_J(Z^h+W)= \\
                                                                \\
(d\>\imath_{\pi_1(A)}\pi_2(sB))^h_J(Z^h)+ (\pi_2(WB))_p(\pi_1(A))
\end{array}
$
\vspace{0.2cm}
\item[$(iii)$] \ \\
$
\begin{array}{lll}
(d\>\imath_{\pi_1({\cal J}_{\alpha}A^h)}\pi_2(B^h))_J(Z^h+W)= \\
                                                               \\
(d\>\imath_{\pi_1(sA)}\pi_2(B))^h_J(Z^h)+ (\pi_2(B))_p(\pi_1(WA))
\end{array}
$
\vspace{0.2cm}
\item[$(iv)$] \ \\
$
\begin{array}{lll}
(d\>\imath_{\pi_1({\cal J}_{\alpha}A^h)}\pi_2({\cal J}_{\alpha}B^h))_J(Z^h+W)= \\
                                                                                \\
(d\>\imath_{\pi_1(sA)}\pi_2(sB))^h_J(Z^h)+ (\pi_2(WB))_p(\pi_1(JA))+
(\pi_2(JB))_p(\pi_1(WA)
\end{array}
$

\end{enumerate}
\end{lemma}

\vspace{0.1cm}

  For any (local) section $a$ of $A(M)$, following \cite{Gaud}, denote by
$\widetilde a$ the vertical vector field on ${\cal G}$ defined by
\begin{equation}\label{eq tilde a}
\widetilde a_J=a+J\circ a\circ J.
\end{equation}
Let us note that for every $J\in {\cal G}$ we can find sections $a_1,...,a_m$,
$m=4n^2-2n$, of $A(M)$ near the point $p=\pi(J)$ such that $\widetilde
a_1,...,\widetilde a_m$ form a basis of the vertical vector space at each point in a
neighbourhood of $J$.

\vspace{0.1cm}

\begin{lemma}\label {H-V brackets}
Let $J\in {\cal G}$ and let $a$ be a section of $A(M)$ near the point $p=\pi(J)$. Then,
for any section $A$ of the bundle $TM\oplus T^{\ast}M$ near $p$, we have:
\begin{enumerate}
\item[$(i)$] $[\pi_1(A^h),\widetilde a]_J=(\widetilde{\nabla_{\pi_1(A)}a})_J$
\vspace{0.2cm}
\item[$(ii)$] $[\pi_1(A^h),{\cal J}_{\alpha}\widetilde a]_J=
(-1)^{\alpha+1}J\circ (\widetilde{\nabla_{\pi_1(A)}a})_J$ \vspace{0.2cm}
\item[$(iii)$] $[\pi_1({\cal J}_{\alpha}A^h),\widetilde
a]_J=(\widetilde{\nabla_{\pi_1(JA)}a})_J-\pi_1(\widetilde a(A))^h_J$
\vspace{0.2cm}
\item[$(iv)$] $[\pi_1({\cal J}_{\alpha}A^h),{\cal J}_{\alpha}\widetilde
a]_J=(-1)^{\alpha+1}[J\circ (\widetilde{\nabla_{\pi_1(JA)}a})_J-\pi_1((J\circ\widetilde
a)(A))^h_J]$
\end{enumerate}
\end{lemma}

\begin{proof}
Let $a(Q_i)=\sum_{j=1}^{4n}\varepsilon_ia_{ij}Q_j$, $i=1,...,4n$. Then, in the local
coordinates of $A(M)$ introduced above,
$$
\widetilde a=\sum_{i<j}\widetilde a_{ij}\frac{\partial}{\partial y_{ij}}
$$
where
$$
\widetilde
a_{ij}=a_{ij}\circ\pi+\sum_{k,l=1}^{4n}y_{ik}(a_{kl}\circ\pi)y_{lj}\varepsilon_k\varepsilon_l.
$$
Let us also note that for every vector field $X$ on $M$ near the point $p$, in view of
(\ref{eq 3.1}), we have
$$
[X^h,\frac{\partial}{\partial y_{ij}}]_J=0,~X_J^h=\sum_{i=1}^{2n}
X^i(p)\frac{\partial}{\partial\tilde x_i}(J),
~(\nabla_{X_p}a)(Q_i)=\sum_{j=1}^{4n}\varepsilon_iX_p(a_{ij})Q_j
$$
since $\nabla Q_i|_p=0$ and $\nabla S_{ij}|_p=0$. Now the lemma follows by simple
computations making use of (\ref{cal J/ver}) and (\ref{cal J/hor}).
\end{proof}

\vspace{0.1cm}

The proof of the next lemma is easy and will be omitted.

\begin{lemma}\label {H-V Lie der&diff}
Let $A$ be a section of the bundle $TM\oplus T^{\ast}M$ and $V$ a vertical vector field
on ${\cal G}$. Then:
\begin{enumerate}
\item[$(i)$] ${\cal L}_{V}\pi_2(A^h)=0$; \hspace{0.2cm}$\imath_{V}\pi_2(A^h)=0$.
\vspace{0.2cm}
\item[$(ii)$] ${\cal L}_{V}\pi_2({\cal J}_{\alpha}A^h)=\pi_2(VA)^h$;
\hspace{0.2cm} $\imath_{V}\pi_2({\cal J}_{\alpha}A^h)=0$.
\end{enumerate}
\end{lemma}

\vspace{0.1cm}

{\it Notations}. Take a point $J\in {\cal G}$ and fix a basis $\{U_{2r-1},U_{2r}={\cal
J}_1U_{2r-1}\}$, $r=1,...,2n^2-n$, of the vertical space ${\cal V}_J$. Now let us take
sections $a_{2r-1}$ of $A(M)$ near the point $p=\pi(J)$ such that $a_{2r-1}=U_{2r-1}$
and $\nabla a_{2r-1}|_p=0$. Define vertical vector fields $\widetilde a_{2r-1}$ by
(\ref{eq tilde a}). Then $\{\widetilde a_{2r-1},{\cal J}_1\widetilde a_{2r-1}\}$,
$r=1,...,2n^2-n$, is a frame of the vertical bundle on ${\cal G}$ near the point $J$.
Denote by $\{\beta_{2r-1},\beta_{2r}\}$ the dual frame. Then $\beta_{2r}={\cal
J}_1\beta_{2r-1}$.

  Under these notations, identity (\ref{eq 3.2}) and Lemmas~\ref{brackets} and \ref{H-V
brackets} imply the following

\vspace{0.1cm}
\begin{lemma}\label {H-Vstar Lie der}
Let $A$ be a section of the bundle $TM\oplus T^{\ast}M$ near the point $p=\pi(J)$. Then
for every $Z\in T_pM$ and $s,q=1,....,4n^2-2n$, we have
\begin{enumerate}
\item[$(i)$]$({\cal L}_{\pi_1(A^h)}\beta_s)_J(Z^h+U_q)=-\beta_s(R(\pi_1(A),Z)J).$
\vspace{0.2cm}
\item[$(ii)$]$({\cal L}_{\pi_1({\cal
J}_{\alpha}A^h)}\beta_s)_J(Z^h+U_q)=-\beta_s(R(\pi_1(JA),Z)J).$
\end{enumerate}
\end{lemma}

\vspace{0.1cm}
\begin{prop}\label {Nijenhuis}
Suppose that the connection $\nabla$ is torsion-free and let $J\in {\cal G}$, $A,B\in
T_{\pi(J)}M\oplus T_{\pi(J)}^{\ast}M$, $V,W\in {\cal V}_J$, $\varphi,\psi\in {\cal
V}_J^{\ast}$. Then:
\begin{enumerate}
\item[$(i)$] \ \\
$
\begin{array}{lll}
N_{\alpha}(A^h,B^h)_J &=&-R(\pi_1(A),\pi_1(B))J+(-1)^{\alpha}J\circ R(\pi_1(A),\pi_1(JB))J\\
                                                                  \\
& &(-1)^{\alpha}J\circ R(\pi_1(JA),\pi_1(B))J+R(\pi_1(JA),\pi_1(JB))J
                                                                   \\
& &\hspace{3.5cm} -\displaystyle{\frac{1}{2}[1+(-1)^{\alpha}]} \omega_{A,B},
\end{array}
$
where $\omega_{A,B}$ is the vertical $1$-form on ${\cal G}$ given by
$$
\begin{array}{lll}
\omega_{A,B}(W)&=&\pi_2(JA)(\pi_1(WB))+\pi_2(WB)(\pi_1(JA)) \\
               & &-\pi_2(JB)(\pi_1(WA))-\pi_2(WA)(\pi_1(JB))
\end{array}
$$
for every $W\in {\cal V}_J$.
\item[$(ii)$] \ \\
$N_{\alpha}(A^h,V)_J=[1+(-1)^{\alpha}]((J\circ V)A)^h_J$

\item[$(iii)$] \ \\
$N_{\alpha}(A^h,\varphi)_J\in {\cal H}_J\oplus {\cal H}_J^{\ast}$ ~ and \ \\
              \ \\
$
\begin{array}{lll}
<\pi_{\ast}N_{\alpha}(A^h,\varphi)_J,B>=-\frac{1}{2}\varphi(N_{1}(A^h,B^h)_J)\\
                                                                              \\
 -\frac{1}{2}[1+(-1)^{\alpha}]
\varphi(J\circ R(\pi_1(A),\pi_1(JB))J+J\circ R(\pi_1(JA),\pi_1(B))J).
\end{array}
$
\item[$(iv)$] \ \\
$N_{\alpha}(V+\varphi,W+\psi)_J=0$.
\end{enumerate}
\end{prop}
\begin{proof}
 Set $p=\pi(J)$ and extend the vectors $A,B$ to (local) sections of $TM\oplus T^{\ast}M$,
denoted again by $A,B$, in such a way that $\nabla A|_p=\nabla B|_p=0$.

 Let $s$ be the section of ${\cal G}$ defined above with the property that
$s(p)=J$ and $\nabla s|_p=0$  ($s$ being considered as a section of $A(M)$).

 According to Lemmas~\ref{brackets}, \ref{Lie deriv} and \ref{Half-diff}, the
part of $N_{\alpha}(A^h,B^h)_J$ lying in  ${\cal H}_J\oplus {\cal H}^{\ast}_J$ is given
by
$$
({\cal H}\oplus {\cal
H}^{\ast})N_{\alpha}(A^h,B^h)_J=(-[A,B]-s[A,sB]-s[sA,B]+[sA,sB])^h_J.
$$
Note that we have $\nabla\pi_1(A)|_p=\pi_1(\nabla A|_p)=0$  and  $\nabla\pi_1(sA)|_p=
\pi_1((\nabla s)|_p(A)+ s(\nabla A|_p))=0$. Similarly, $\nabla\pi_2(A)|_p=0$ and
$\nabla\pi_2(sA)|_p=0$. We also have $\nabla\pi_1(B)|_p=0$, $\nabla\pi_1(sB)|_p=0$ and
$\nabla\pi_2(B)|_p=0$, $\nabla\pi_2(sB)|_p=0$. Now, since $\nabla$ is torsion-free, we
easily get $({\cal H}\oplus {\cal H}^{\ast})N_{\alpha}(A^h,B^h)_J=0$ by means of the
following simple observation: Let $Z$ be a vector field and $\omega$ a $1$-form on $M$
such that $\nabla Z|_p=0$ and $\nabla\omega|_p=0$. Then for every $T\in T_pM$
$$
({\cal L}_{Z}\omega)(T)_p =(\nabla_{Z}\omega)(T)_p=0 ~\mbox { and
}~(d\,\imath_{Z}\omega)(T)_p=(\nabla_{T}\omega)(Z)_p=0.
$$

  By Lemmas~\ref{brackets}, \ref{Lie deriv} and \ref{Half-diff} the vertical part of
$N_{\alpha}(A^h,B^h)_J$ is equal to
$$
\begin{array}{lll}
{\cal V}N_{\alpha}(A^h,B^h)_J&=&
  -R(\pi_1(A),\pi_1(B))J - {\cal J}_{\alpha}R(\pi_1(A),\pi_1(JB))J \\
 & &-{\cal J}_{\alpha}R(\pi_1(JA),\pi_1(B))+R(\pi_1(JA),\pi_1(JB))J\\
\end{array}
$$
The part of $N_{\alpha}(A^h,B^h)_J$ lying in ${\cal V}^{\ast}_J$ is the $1$-form whose
value at every vertical vector $W$ is
$$
\begin{array}{lll}
({\cal V}^{\ast}N_{\alpha}(A^h,B^h)_J)(W)= \\
-\frac{1}{2}[\pi_2(JA)(\pi_1(WB))+\pi_2(WB)(\pi_1(JA)) \\
\hspace{2.5cm}+(-1)^{\alpha}\pi_2(B)(\pi_1((J\circ W)A))+(-1)^{\alpha}\pi_2((J\circ W)A)(\pi_1(B))] \\
                                                                           \\
+\frac{1}{2}[\pi_2(JB)(\pi_1(WA))+\pi_2(WA)(\pi_1(JB)) \\
\hspace{2.5cm}+(-1)^{\alpha}\pi_2(A)(\pi_1((J\circ W)B))+(-1)^{\alpha}\pi_2((J\circ
W)B)(\pi_1(A))].
\end{array}
$$
The endomorphism $W$ of $T_pM\oplus T^{\ast}_pM$ is skew-symmetric with respect to the
natural neutral metric and anti-commutes with $J$, so $<JA,WB>=\\ <(J\circ W)A,B>$.
This gives
$$
\begin{array}{c}
\pi_2(JA)(\pi_1(WB))+\pi_2(WB)(\pi_1(JA))= \\
\pi_2(B)(\pi_1((J\circ W)A))+\pi_2((J\circ W)A)(\pi_1(B)).
\end{array}
$$
It follows that
$$
\begin{array}{lll}
{\cal V}^{\ast}N_{\alpha}(A^h,B^h)_J=-\frac{1}{2}[1+(-1)^{\alpha}]
[\pi_2(JA)(\pi_1(WB))+\pi_2(WB)(\pi_1(JA)) \\
\hspace{5.8cm}               -\pi_2(JB)(\pi_1(WA))-\pi_2(WA)(\pi_1(JB))].
\end{array}
$$
This proves $(i)$.

  To prove $(ii)$ take a section $a$ of $A(M)$ near the point $p$ such that $a(p)=V$
and $\nabla a|_p=0$. Let $\widetilde a$ be the vertical vector field defined by
(\ref{eq tilde a}). Then it follows from Lemmas~\ref{H-V brackets} and \ref{H-V Lie
der&diff} that
$$
N_{\alpha}(A^h,V)_J=\frac{1}{2}N_{\alpha}(A^h,\widetilde a)_J= ((J\circ
V)A+(-1)^{\alpha}(J\circ V)A)^{h}_J.
$$

  To prove $(iii)$ let us take the vertical coframe $\{\beta_{2r-1},\beta_{2r}\}$,
$r=1,...,2n^2-n$, defined before the statement of Lemma~\ref{H-Vstar Lie der}. Set
$\varphi=\sum_{s=1}^{4n^2-2n}\varphi_s\beta_s$, $\varphi_s\in {\Bbb R}$. Let $E_1,...,
E_{2n}$  be a basis of $T_pM$ and $\xi_1,...,\xi_{2n}$ its dual basis. Then, by
Lemma~\ref{H-Vstar Lie der}, we have \vspace{0.1cm}
\begin{equation}\label{eq H-Vstar Nij}
\begin{array}{lll}
N_{\alpha}(A^h,\varphi)_J=\sum_{s=1}^{4n^2-2n}\varphi_s N_{\alpha}(A^h,\beta_s)_J=\\
                                                                                   \\
\sum_{s=1}^{4n^2-2n}\sum_{k=1}^{2n}\varphi_s\{[\beta_s(R(\pi_1(A),E_k)J)+
(-1)^{\alpha+1}\beta_s(J\circ R(\pi_1(JA),E_k)J)](\xi_k)^h_J \\
                                                                                    \\
+[\beta_s(R(\pi_1(JA),E_k)J)-(-1)^{\alpha+1}\beta_s(J\circ
R(\pi_1(A),E_k)J)](J\xi_k)^h_J\}.
\end{array}
\end{equation}
Moreover, note that
$$
<\xi_k,B>=\frac{1}{2}\xi_k(\pi_1(B)) \mbox { and }
<J\xi_k,B>=-<\xi_k,JB>=-\frac{1}{2}\xi_k(\pi_1(JB)),
$$
therefore
$$
\sum_{k=1}^{2n}<\xi_k,B>E_k=\frac{1}{2}\pi_1(B) \mbox { and }
\sum_{k=1}^{2n}<J\xi_k,B>E_k=-\frac{1}{2}\pi_1(JB).
$$
Now $(iii)$ is an obvious consequence of (\ref{eq H-Vstar Nij}) and formula $(i)$.

  Finally, identity $(iv)$ follows from the fact that the generalized almost complex structure
${\cal J}_{\alpha}$ on every fibre of ${\cal G}$ is induced by a complex structure.
\end{proof}

\section{The integrability condition}

\begin{tw}\label {Integrability}
Let $M$ be a $2n$-dimensional manifold and $\nabla$ a torsion-free linear connection on
$M$. Let ${\cal J}_{\alpha}={\cal J}_{\alpha}^{\nabla}$, $\alpha=1,2$, be the
generalized almost complex structures induced by $\nabla$ on the twistor space ${\cal
G}$ of $M$. Then:
\begin{enumerate}
\item[$(i)$] If $n=1$, the  structure ${\cal J}_1$ is always
integrable.
\item[$(ii)$] If $n\geq 2$, the structure ${\cal J}_1$ is integrable if and only if the
connection $\nabla$ is flat.
\item[$(iii)$] The structure ${\cal J}_2$ is never integrable.
\end{enumerate}
\end{tw}

\begin{proof}
 $(i)$~Let $J\in {\cal G}$ and let $\{Q_1, Q_2=JQ_1, Q_3, Q_4=JQ_3\}$ be an orthonormal basis
of $T_pM\oplus T_p^{\ast}M$, $p=\pi(J)$. By Proposition~\ref{Nijenhuis}, the Nijenhuis
tensor $N_1$ of ${\cal J}_1$ vanishes at the point $J$ if and only if
$N_1(Q_1^h,Q_3^h)=0$ (in view of the fact that $N_1({\cal J}_1E,F)=N_1(E,{\cal
J}_1F)=-{\cal J}_1N(E,F)$ for $E,F\in T{\cal G}$).

   Let $\pi_1(Q_i)=e_i$, $i=1,...,4$. Then, according to Proposition~\ref{Nijenhuis},
$$
N_1(Q_1^h,Q_3^h)=([-R(e_1,e_3)J+R(e_2,e_4)J]-J\circ[R(e_1,e_4)J+R(e_2,e_3)J])^h_J.
$$
Both summands in the above formula vanishes since, by Lemma~\ref{Orient},
$e_3=\cos{t}\,e_1+\sin{t}\,e_2$, $e_4=-\sin{t}\, e_1+\cos{t}\,e_2$ for some $t\in {\Bbb
R}$.

\vspace{0.1cm}

$(ii)$~ Let $dim\, M=2n\geq 4$ and assume that the generalized almost complex structure
${\cal J}^{\nabla}_1$ is integrable. Then, by Proposition~\ref{Nijenhuis}, for every
$p\in M$, every (genuine) complex structure $K$ on $T_pM$ and every $X,Y\in T_pM$ we
have
$$
R(X,Y)K+K\circ R(X,KY)K+K\circ R(KX,Y)K-R(KX,KY)K=0,
$$
where $R$ is the curvature tensor of the connection $\nabla$. The latter identity, as
is well-known, is the integrability condition for the Atiyah-Hitchin-Singer almost
complex structure on the twistor space of complex structures on the tangent spaces of
$M$ (see, e.g., \cite[Theorem 1 or Theorem 3]{OBrRaw}). Then, according to the
arguments of \cite[pp. 42 - 43]{OBrRaw}, there exists a bilinear form $\mu$ on
$TM\times TM$ such that
\begin{equation}\label{eq curv}
R(X,Y)Z= \mu(X,Y)Z-\mu(Y,X)Z+\mu(X,Z)Y-\mu(Y,Z)X.
\end{equation}
Now let $p\in M$ and let $\{E_1,...,E_{2n}\}$ be an arbitrary basis of $T_pM$. Denote
by $\{\eta_1,...,\eta_{2n}\}$ its dual basis. Let $J$ be the complex structure on
$T_pM\oplus T^{\ast}_pM$ for which $JE_{2k-1}=\eta_{2k}$, $JE_{2k}=-\eta_{2k-1}$,
$k=1,...,n$. This structure is compatible with the natural neutral metric of
$T_pM\oplus T^{\ast}_pM$, so we get a generalized complex structure sending $T_pM$ onto
$T^{\ast}_pM$ and vice versa (it is similar to the structure in Example 2, Section 2).
If $n$ is an even number, the structure $J$ yields the canonical orientation of
$T_pM\oplus T^{\ast}_pM$, hence $J\in {\cal G}$. So, suppose that $n$ is even and let
$X,Y\in T_pM$. Then, by Proposition~\ref{Nijenhuis}\,$(i)$, we have $R(X,Y)J=0$. The
latter identity is equivalent to the following identities
$$
R(X,Y)\eta_{2k}-JR(X,Y)E_{2k-1}=0,~ R(X,Y)\eta_{2k-1}+JR(X,Y)E_{2k}=0,
$$
$k=1,...,n$. It follows that for every $Z\in T_pM$, we have
$$
2[\mu(X,Y)-\mu(Y,X)]\eta_{2k}(Z)+\mu(X,Z)\eta_{2k}(Y)-\mu(Y,Z)\eta_{2k}(X) $$
$$
+\mu(X,E_{2k-1})(JY)(Z)-\mu(Y,E_{2k-1})(JX)(Z)=0
$$
and
$$
2[\mu(X,Y)-\mu(Y,X)]\eta_{2k-1}(Z)+\mu(X,Z)\eta_{2k-1}(Y)-\mu(Y,Z)\eta_{2k-1}(X) $$
$$
-\mu(X,E_{2k})(JY)(Z)+\mu(Y,E_{2k})(JX)(Z)=0.
$$
Let $k\neq l$, $1\leq k,l\leq n$, be two indexes. Putting $X=E_{2k-1}$, $Y=E_{2k}$,
$Z=E_{2l-1}$ into the above identities, we get $\mu(E_{2k-1},E_{2l-1})=0$ and
$\mu(E_{2k},E_{2l-1})=0$; putting $X=E_{2k}$, $Y=E_{2k-1}$, $Z=E_{2l}$, we get
$\mu(E_{2k-1},E_{2l})=0$ and $\mu(E_{2k},E_{2l})=0$. Moreover, setting $X=Z=E_{2l-1}$,
$Y=E_{2k}$ and $X=Z=E_{2l}$, $Y=E_{2k}$ into the first of the above identities, we
obtain $\mu(E_{2l-1},E_{2l-1})=0$ and $\mu(E_{2l},E_{2l})=0$. Now take the basis
$E_1'=E_1$, $E_2'=E_3$, $E_3'=E_2$, $E_4'=E_4$,...,$E_{2n}'=E_{2n}$. Then the
identities $\mu(E_1',E_3')=\mu(E_3',E_1')=0$ give $\mu(E_1,E_2)=\mu(E_2,E_1)=0$. It
follows that $\mu(E_{2k-1},E_{2k})=\mu(E_{2k},E_{2k-1})=0$. Therefore $\mu=0$, thus
$R=0$ by (\ref{eq curv}).

  Now assume that $n=2m+1$ is an odd number. Let $X,Y\in T_pM$ be two linearly
independent tangent vectors. Let  $\{E_1,...,E_{2n}\}$ be an arbitrary basis of $T_pM$
with $E_1=X$, $E_2=Y$. Denote by $\{\eta_1,...,\eta_{2n}\}$ its dual basis. Let $J$ be
the complex structure on $T_pM\oplus T^{\ast}_pM$ for which $JE_{2k-1}=\eta_{2k}$,
$JE_{2k}=-\eta_{2k-1}$, $k=1,...,m$, and $JE_{4m+1}=E_{4m+2}$,
$J\eta_{4m+1}=\eta_{4m+2}$. Then $J\in {\cal G}$ and the preceding arguments show that
$$
R(X,Y)E_i=0 \mbox { for } i=1,...,4m.
$$
Applying the latter identity for the basis $\{E_1, E_2, E_{4m+1}, E_{4m+2},
E_3,...,E_{4m}\}$, we see that $R(X,Y)E_{4m+1}=R(X,Y)E_{4m+2}=0$. It follows that
$R=0$.

\vspace{0.1cm}

$(iii)$~ Let $J\in {\cal J}$ and let $\{Q_1, Q_2=JQ_1,..., Q_{4n-1},
Q_{4n}=JQ_{4n-1}\}$ be an orthonormal basis of $T_pM\oplus T_p^{\ast}M$, $p=\pi(J)$.
Define endomorphisms $S_{ij}$ of $T_pM\oplus T_p^{\ast}M$ by (\ref{eq Sij}). Then
$V=S_{13}-S_{24}$ is a vertical vector of ${\cal G}$ at the point $J$ and it follows
from Proposition~\ref{Nijenhuis}$\,(ii)$ that $N_{2}(E_1^h,V)_J=2(Q_4)_J^h\neq 0$.
\end{proof}

\vspace{0.1cm}

\noindent {\bf Remarks.} {\bf 1}. Concerning Theorem 1$(ii)$, let us note that that
every flat torsion-free connection on a manifold induces an affine structure on it,
i.e. local coordinates whose transition functions are affine. If, in addition, the
connection is complete, then the manifold is the quotient of an affine space  by a
group of affine transformations acting freely and properly discontinuously on it (see,
for example, \cite{Wolf}).

\noindent {\bf 2}. The complex structure ${\cal K}$ on the fibres of ${\cal G}$ is
K\"ahlerian with respect to the metric $<a,b>=\displaystyle{-\frac{1}{2}}Trace\,(a\circ
b)$. Let ${\cal S}$ be the generalized complex structure on the vertical spaces of
${\cal G}$ induced by the K\"ahler form of K\"ahler structure $({\cal J},<~,~>)$ (see
Example 2, Section 2). Then, given a connection $\nabla$ on $M$, we can define two new
generalized almost complex structures ${\cal I}_{\alpha}$ on the twistor space ${\cal
G}$ by setting ${\cal I}_{\alpha}=(-1)^{\alpha+1}{\cal S}$ on ${\cal V}\oplus {\cal
V}^{\ast}$ and ${\cal I}_{\alpha}={\cal J}_{1}^{\nabla}\,(={\cal J}_{2}^{\nabla})$ on
${\cal H}\oplus {\cal H}^{\ast}$, $\alpha=1,2$. It is easy to see that the structures
${\cal I}_{\alpha}$ are never integrable. Indeed, let us adopt the notations used in
part $(iii)$ of the proof of Theorem~\ref{Integrability}. Then Lemma~\ref{H-V brackets}
implies that the projection of the Nijenhuis tensor $N_{{\cal I}_{\alpha}}(E_1^h,V)_J$
onto ${\cal H}_J$ is equal to $(Q_4)_J^h$.

\end{document}